\documentclass{amsart}            

\usepackage{amssymb,amsmath}
\usepackage{color}
\usepackage{pdfsync}

\theoremstyle{plain}
\newtheorem{theorem}{Theorem}[section]
\newtheorem*{theorem*}{Theorem}
\newtheorem{lemma}[theorem]{Lemma}
\newtheorem{corollary}[theorem]{Corollary}

\theoremstyle{remark}
\newtheorem{remark}[theorem]{Remark}
\numberwithin{equation}{section}

\newcount\quantno
\everydisplay{\quantno=0}\everycr{\quantno=0}
\newcommand\quant{\advance\quantno by1
                      \ifnum\quantno=1\qquad\else\quad\fi\forall }

\newcommand\rest[1]{\kern-.1em
          \lower.5ex\hbox{$\scriptstyle #1$}\kern.05em}

\newcommand\Bigset[1]{\Bigl\{#1\Bigr\}}

\renewcommand\mod[1]{\left\vert{#1}\right\vert}
\newcommand\bigmod[1]{\bigl\vert{#1}\bigr|}
\newcommand\Bigmod[1]{\Bigl\vert{#1}\Bigr|}

\newcommand\norm[2]{{\Vert{#1}\Vert_{#2}}}

\newcommand\prodo[2]{\left\langle#1,#2\right\rangle}

\newcommand\1{\mathbf{1}}

\newcommand\wrt{\,\text{\rm d}}

\newcommand\BC{\mathbb{C}}

\newcommand\BN{\mathbb{N}}
\newcommand\BR{\mathbb{R}} \newcommand\BRd{\mathbb{R}^d}

\newcommand\cB{\mathcal{B}}

\newcommand\cL{\mathcal{L}}

\newcommand\cP{\mathcal{P}}

\newcommand\ga{\gamma}

\newcommand\la{\lambda}

\newcommand\funnyk{k\hbox to 0pt{\hss\phantom{g}}}

\newcommand\lu[1]{L^1(#1)}
\newcommand\hu[1]{H^1(#1)}
\newcommand\lp[1]{L^p(#1)}

\newcommand\ld[1]{L^2(#1)}

\newcommand\ly[1]{L^\infty(#1)}

\newcommand\BMO[1]{BMO(#1)}

\newcommand\whH{\widehat{\phantom{G}}\hbox to 0pt{\hss $H$}}

\newcommand\emspace{\hbox to 6pt{\hss}}

\newcommand\rmii{\hbox{\rm (ii)}}

\newcommand\OU{{Ornstein--Uhlenbeck}}

\newcommand\e{\mathrm{e}}

\begin{document}

\date{}

\title[Endpoint estimates for Riesz transforms]
{Endpoint estimates \\ for first-order 
Riesz transforms \\ associated 
to the Ornstein--Uhlenbeck operator}

\subjclass[2000]{} 


\thanks{Work partially supported by the
Progetto Cofinanziato ``Analisi Armonica''.}

\author[G. Mauceri,  S. Meda and P. Sj\"ogren]
{Giancarlo Mauceri, Stefano Meda and Peter Sj\"ogren}

\address{Dipartimento di Matematica \\
Universit\`a di Genova \\ via Dodecaneso 35, 16146 Genova \\ Italia}
\email{mauceri@dima.unige.it}
\address{Dipartimento di Matematica e Applicazioni
\\ Universit\`a di Milano-Bicocca\\
via R.~Cozzi 53\\ 20125 Milano\\ Italy}
\email{stefano.meda@unimib.it}
\address{Mathematical Sciences \\
University of Gothenburg and Mathematical Sciences \\Chalmers 
 \\ SE-412 96 G\"oteborg \\ Sweden}
\email{peters@chalmers.se}

\begin{abstract}
In the setting of  Euclidean space with the
Gaussian measure $\ga$, we consider all first-order
Riesz transforms associated to the infinitesimal generator of
the \OU\ semigroup. These operators are known to be 
bounded on $L^p(\ga)$, for
$1<p<\infty$. We determine which of them are bounded from $\hu{\ga}$ to
$\lu{\ga}$ and from $L^\infty(\ga)$ to $\BMO{\ga}$. Here  $\hu{\ga}$ and  $\BMO{\ga}$
are the spaces introduced in this setting by the  first two authors.
Surprisingly, we find that the results depend on the dimension
of the ambient space. 
\end{abstract}

\maketitle

\section{Introduction}

Denote by $\ga$ the normalised Gaussian measure on $\BR^d$ having
density $x\mapsto\pi^{-d/2}\, \e^{-\mod{x}^2}$ 
with respect to the Lebesgue measure. 
The \OU\ operator $\cL$ is the closure in $\ld{\ga}$ of the operator $\cL_0$,
defined on $C^\infty_c(\BR^n)$ by
$$
\cL_0
= -\frac{1}{2}\,  \Delta + x\cdot\nabla,
$$
where $\Delta$ and $\nabla$ denote the Euclidean
Laplacian and gradient, respectively.
The operator $\cL$ is self-adjoint, and its
spectral resolution is
$$
\cL 
= \sum^{\infty}_{j=0} j\, \cP_j,
$$
where $\cP_j$ denotes the orthogonal projection onto the linear span 
of Hermite polynomials of degree~$j$ in $d$ variables. 
The \OU\ operator generates a diffusion semigroup,
which has been the object of many investigations
during the last two decades.  
In particular, efforts have been made
to study operators related to the \OU\ semigroup, with emphasis
on maximal operators \cite{S,GU,MPS, GMMST2},
Riesz transforms 
\cite{Mu, Gun, M, Pisier, Pe, Gut, GST, FGS, FoS, GMST1, PS, U, DV}
and functional calculus \cite{GMST2,GMMST1,MMS1,HMM}.

In this paper, we shall concentrate on the first-order Riesz transforms
associated with $\cL$.  Since the eigenspace associated to the 
zero eigenvalue is nontrivial, 
the positive square root of $\cL$ is not invertible, and one must
define the Riesz transforms carefully.
Consider the sequence $M:\BN\rightarrow \BC$, defined by
$$
M(j)=
\begin{cases}
0       & {\rm if}\ j=0\\
j^{-1/2}  & {\rm if}\ j=1,2,\ldots
\end{cases}
$$
The operator $M(\cL)$, spectrally defined on $\ld{\ga}$,   
extends to a bounded operator on $\lp{\ga}$ for every $p$ in $(1,\infty)$
(see, for instance, \cite[Theorem~1.2]{MMS1}). 
The operator $M(\cL)$ plays, in this setting, the
same role as the operator $(-\Delta)^{-1/2}$
in Euclidean harmonic analysis.  
We denote by $\partial_i$ the differentiation operator with respect to
the variable $x_i$. The formal adjoint of $\partial_i$
in $\ld{\ga}$ is the operator $\partial_i^* = 2x_i - \partial_i$. 
Notice that, at least formally, 
$$
\sum_{i=1}^d \partial_i^* \partial_i = 2\, \cL. 
$$
For each $i$ in $\{1,\ldots,d\}$,
we define the operators $R_i$, $S_i$, $R_i^*$ and $S_i^*$ 
on finite linear combinations of Hermite polynomials by
\begin{equation} \label{f: def Riesz}
\begin{array}{ll}
R_i 
 = \partial_i M(\cL)  
& \qquad S_i 
 =  M(\cL)  \partial_i  \\
R_i^* 
 = M(\cL) \partial_i^* 
& \qquad S_i^* 
=  \partial_i^* M(\cL). 
\end{array}
\end{equation}
It is straightforward to check that $R_i^*$ and $S_i^*$
are the \emph{formal} adjoints in $\ld{\ga}$ of $R_i$ and $S_i$,
respectively.  
Recall the action of $\partial_i$ and $\partial_i^*$
on Hermite polynomials:
$$
\partial_i H_n 
= 2n\, H_{n-1}
\qquad\hbox{and}\qquad
\partial_i^* H_n 
= H_{n+1}, 
$$
where $H_j$ denotes the Hermite polynomials of degree $j$
in the $i^{\small\textrm{th}}$ variable.  
A straightforward argument, which uses these formulae and 
the expression for the $\ld{\ga}$ norm of a Hermite polynomial, 
then shows that $R_i$, $R_i^*$, $S_i$ and $S_i^*$ extend to
bounded operators on $\ld{\ga}$, and that $R_i^*$ and $S_i^*$
are the Hilbert space adjoints of $R_i$ and $S_i$,
respectively.   
Note that, in contrast with the Euclidean case, these transforms are
not antisymmetric. 

The operators $R_i$, $S_i$, $R_i^*$
and $S_i^*$ are all bounded on $\lp{\ga}$ for each $p$ in $(1,\infty)$
(see \cite{M}).
This may be easily proved by using the commutation relations
between $\partial_i$ and $\cL$ and the spectral multiplier
result \cite[Theorem~1.2]{MMS1}.
It is straightforward to check that none of these operators 
are bounded on $\lu{\ga}$ or on $\ly{\ga}$.  
Weak type $1$ estimates for $R_i$ and for $S_i^*$ have been proved in \cite{FGS}
and \cite{AFS}, respectively.

Thus, it is natural to ask what further estimates
these operators satisfy at the endpoints $p=1$ and $p=\infty$.
In the Euclidean case there are 
substitute results at $p=1$ and $p=\infty$ saying that
the operators  are bounded from $\hu{\BR^d}$ to $\lu{\BR^d}$ and
from $\ly{\BR^d}$ to $\BMO{\BR^d}$.
In our setting, the  spaces  
$\hu{\ga}$ and $\BMO{\ga}$ were defined by Mauceri and Meda
in \cite{MM}, where these authors also developed a theory of singular 
integral operators in the Gaussian setting (see also \cite{CMM}
for a related theory in a more general setting).  As 
applications, they proved that the imaginary
powers of the \OU\ operator are bounded from $\hu{\ga}$ 
to $\lu{\ga}$ and from $\ly{\ga}$ to $\BMO{\ga}$, and that 
the operators $R_i$ are bounded from $\ly{\ga}$ to $\BMO{\ga}$.

The question to be studied below is whether the
Riesz transforms defined in (\ref{f: def Riesz})
are bounded from  $\hu{\ga}$ to $\lu{\ga}$ and
from $\ly{\ga}$ to $\BMO{\ga}$. 
We find that this boundedness always holds in dimension one. 
But, surprisingly, in higher dimensions, each Riesz
transform has one of these boundedness properties but not the other.
Our results for $d\geq 2$ are summarised in the table
in Theorem~\ref{t: main}.

The paper is organised as follows.  Section 2 contains
background material and a few preliminary results,
including the expressions for the kernels of $R_i$ and $S_i^*$. 
Our main result, Theorem~\ref{t: main}, is stated 
in Section~\ref{s: TMR}.  Its proof constitutes Sections~\ref{R}-\ref{D}.  

\section{Preliminaries}

We briefly recall the Hardy space $\hu{\ga}$ from \cite{MM}.
A Euclidean ball $B$ is called \emph{admissible} if 
\begin{equation} \label{f: 1-admissible} 
r_B \leq \min\bigl(1,1/\mod{c_B}\bigr);
\end{equation}
here and in the sequel
$r_B$ and $c_B$ denote the radius and the centre of $B$, respectively. 
The collection of all admissible balls will be denoted by $\cB_1$.
When there is equality in (\ref{f: 1-admissible}), we say that $B$
is a \emph{maximal ball} in $\cB_1$. 

A (Gaussian) \emph{atom} is either the constant function $1$ or 
a function $a$ in $\ly{\ga}$,
supported in an admissible ball $B$ and such that
\begin{equation}  \label{f: prop atom}
\norm{a}{\infty} \leq \ga(B)^{-1}
\qquad\hbox{and}\qquad
\int_{\BR^n} a \, \wrt \ga =0.
\end{equation}
The space $\hu{\ga}$ is then the vector space of all functions 
$f$ in $\lu{\ga}$ that admit a decomposition of the form 
$\sum_j \la_j \, a_j$, where the $a_j$ are atoms 
and the sequence of complex numbers $\{\la_j\}$ is summable. 
The norm of $f$ in $\hu{\ga}$ is defined as the infimum
of $\sum_j\mod{\la_j}$ over all representations of $f$ as above.
In \cite{MM} the space $\hu{\ga}$ is defined by means of $(1,p)$-atoms
with $1<p<\infty$, but in \cite{MMS2} it is verified
that the space obtained is the same as ours.

Note that $\hu{\ga}$ is defined much as the atomic space $H^1$
on spaces of homogeneous type in the sense
of R.R.~Coifman and G.~Weiss \cite{CW}, but with one difference. 
Namely, 
only the exceptional atom and atoms with ``small supports'',
i.e., with supports contained in admissible balls,
appear in the definition of $\hu{\ga}$.
This difference is quite significant and has important consequences.   

It is known \cite[Theorem~5.2]{MM} 
that the Banach dual of $\hu{\ga}$ is isomorphic
to the space $BMO(\ga)$ of all functions of ``bounded mean 
oscillation'', i.e., of all functions in $\lu{\ga}$ such that 
\begin{equation} \label{f: BMO}
\norm{f}{*}
= \sup_{B \in \cB_1}  \frac{1}{\ga(B)} \int_B \mod{f-f_B} \wrt \ga
< \infty,
\end{equation}
where $f_B = \frac{1}{\ga(B)} \int_B f \wrt \ga$.
A convenient norm on $BMO(\ga)$ is the following
$$
\norm{f}{BMO}
= \norm{f}{1} + \norm{f}{*}.
$$

\begin{remark} \label{rem: eq norm BMO}
An equivalent norm on $BMO(\ga)$ is obtained by replacing in (\ref{f: BMO})
balls in $\cB_1$ by cubes of sidelength at most $2\, \min(1,1/\mod{c_B})$ 
(see \cite[Section~2]{MM}).  
As a consequence, a function $f$ is in $\hu{\ga}$
if and only if it admits a decomposition of the form $\sum_j \la_j\, b_j$
where $\{\la_j\}$ is summable and the $b_j$ are either the constant
function $1$ or functions in $\ly{\ga}$
supported in cubes $Q$ of sidelength at most $2\, \min(1,1/\mod{c_B})$, 
satisfying
$$
\norm{b}{\infty} \leq \ga(Q)^{-1}
\qquad\hbox{and}\qquad
\int_{\BR^n} b \, \wrt \ga =0.
$$
An equivalent norm on
$\hu{\ga}$ is then obtained by taking the infimum 
of $\sum_j\mod{\la_j}$ over all representations of $f$ as above.
\end{remark}

Given a bounded operator $T$ on $\ld{\ga}$, 
we denote by $K_T$ the Schwartz kernel of $T$ and by $k_T$
the kernel of $T$ with respect to the Gaussian measure,
defined by
\begin{equation} \label{f: kT and KT}
k_T(x,y)
= \pi^{d/2} \, \e^{\mod{y}^2}\, K_T(x,y),
\end{equation}
in the sense of distributions in $\BR^d\times \BR^d$. 
The reason for introducing $k_T$ is that if $K_T$ is locally integrable
in $\BR^d\times \BR^d$, then
$T$ is an integral operator with kernel $k_T$ with
respect to $\ga$, i.e.,
$$
Tf(x) 
= \int_{\BR^d} k_T(x,y) \, f(y) \wrt \ga(y)
\quant x \in \BR^d \quant f \in C_c(\BR^d).
$$
One of the results in \cite{MM} says that if $T$
is bounded on $\ld{\ga}$ and $k_T$ is locally integrable 
off the diagonal of $\BR^d \times \BR^d$, and satisfies
\begin{equation} \label{f: HIC H1L1}
\sup_{B\in \cB_1} \, r_B \, \sup_{y\in B} \int_{(2B)^c}
\bigmod{\nabla_y k_T(x,y)} \wrt\ga(x)
< \infty,
\end{equation}
then $T$ is bounded from $\hu{\ga}$ to $\lu{\ga}$,
and, consequently, on $\lp{\ga}$ for all $p$ in $(1,2)$.  
Notice that (\ref{f: HIC H1L1}) 
is a local H\"ormander integral condition.   
Similarly, if
\begin{equation} \label{f: HIC LIBMO}
\sup_{B\in \cB_1} \, r_B \, \sup_{x\in B} \int_{(2B)^c}
\bigmod{\nabla_x k_T(x,y)} \wrt\ga(y)
< \infty,
\end{equation}
then $T$ is bounded from $\ly{\ga}$ to $\BMO{\ga}$,
and, consequently, on $\lp{\ga}$ for all $p$ in $(2,\infty)$.  

We shall use these criteria to prove the 
boundedness results contained in Theorem~\ref{t: main},
and start by determining the kernels of $R_i$
and of $S_i^*$.  The functions $\phi$ and $\psi$, defined~by
\begin{equation*} \label{fipsi}
\phi(r,x,y)
= \frac{ry-x}{\sqrt{1-r^2}} \qquad \mathrm{and} \qquad
\psi(r,x,y)
= \frac{rx-y}{\sqrt{1-r^2}}
\quant x,\,y \in  \BRd \quant r \in (0,1),
\end{equation*}
will occur frequently.
Note that if $d\geq 2$, then $\phi$ and $\psi$ are vector-valued;
their components will be denoted by $\phi_i$ and $\psi_i$.
The arguments of these two functions will often be suppressed.

As shown in \cite[Lemma 2.2]{GMST1}, 
the Schwartz kernel of the operator $M(\cL)$, defined
in the introduction, is 
\[
K_{M(\cL)}(x,y) 
= \pi^{-(d+1)/2} \int_0^1  \Bigl[\frac{\e^{-|\psi|^2}}{(1-r^2)^{d/2}}
 - \e^{-|y|^2}\Bigr]\wrt \rho(r),
\]
where the measure $\rho$, supported in $[0,1]$, is defined by
\begin{equation}  \label{f: measure rho}
\wrt \rho(r) 
= \frac{\wrt r}{r\, \sqrt{-\log r}}.
\end{equation}
Hence the Schwartz kernel $K_{R_i}$ of the operator $R_i$ 
agrees off the diagonal with the function
\[  
\partial_{x_i} K_{M(\cL)}(x,y) 
= - 2\pi^{-(d+1)/2} \int_0^1 
      \frac{r\, \psi_i\, \e^{-|\psi|^2}}{(1-r^2)^{(d+1)/2}} \wrt \rho(r).
\]
By using (\ref{f: kT and KT}) and the fact that
${|\phi|^2-|\psi|^2} =|x|^2-|y|^2 $, we get that 
\begin{equation} \label{eq:kern}
k_{R_i}(x,y)
= -\frac2{\sqrt{\pi}} \, \e^{|x|^2}  \int_0^1 
\frac{r\, \psi_i\,  \e^{-|\phi|^2}}{(1-r^2)^{(d+1)/2}} \wrt \rho(r)
\end{equation}
for all $x\neq y$. 
To determine $k_{S_i^*}$, we argue similarly, observing that 
\[
 \partial_{x_i}^* \e^{-|\psi|^2}
= \frac{-2\phi_i}{\sqrt{1-r^2}} \,\, \e^{-|\psi|^2},
\]
and obtain
\begin{align} \label{eq:kerntilde}
k_{S_i^*}(x,y)
& = -\frac2{\sqrt{\pi}} \, \e^{|y|^2} \int_0^1 
      \Bigl[\frac{\phi_i \, \e^{-|\psi|^2} }{(1-r^2)^{(d+1)/2}} 
      + x_i \, \e^{-|y|^2}\Bigr] \wrt \rho(r) \notag \\
& = -\frac2{\sqrt{\pi}} \, \e^{|x|^2} \int_0^1
      \Bigl[\frac{\phi_i \, \e^{-|\phi|^2} }{(1-r^2)^{(d+1)/2}} 
      + x_i \,  \e^{-|x|^2}\Bigr] \wrt \rho(r)
\end{align}
for all $x\neq y$.

The constants $C<\infty$ and $c>0$ will depend only
on the dimension $d$, and they may vary from occurrence to 
occurrence. By  $f\sim g$ we mean $c < f/g < C$.
Lebesgue measure is denoted   $\la$.

\section{Statements of results}  \label{s: TMR}

The main result of this paper is the following.

\begin{theorem}\label{t: main}
In dimension one, the four operators $R_1$, $R_1^*$, $S_1$ and $S_1^*$ 
are bounded from $\hu{\ga}$ to $\lu{\ga}$ and from $\ly{\ga}$ to $BMO(\ga)$.

In dimension $d\ge2$, the boundedness properties of the operators 
$R_i,R^*_i,S_i$ and $S^*_i$ are given by the following table,
where $B$ means ``bounded" and $U$ ``unbounded".
\medskip
\begin{center}
\begin{tabular}{|c||c|c|c|c|} 
\hline
${\phantom a}$ & $R_i$ & $S_i$&$R_i^*$&$S_i^*$\\
\hline \hline
$H^1 \to L^1$      & $U$ & $U$ & $B$ & $B$ \\
\hline
$L^\infty \to BMO$ & $B$ & $B$ & $U$ & $U$ \\
\hline
\end{tabular}
\end{center}
\medskip
\end{theorem}

\begin{remark}
Note that all the operators $R_i$ and $S_i$, $i$ in $\{1,\ldots,d\}$,
have the same boundedness properties, and so do
the operators $R_i^*$ and $S_i^*$.  Furthermore,
the boundedness properties of $R_i^*$ and $S_i^*$
are ``dual'' to those of $R_i$ and $S_i$.  This is not
an accident, and the proofs of the statements concerning
the boundedness properties of
$R_i^*$ and $S_i$ will be obtained, by a duality argument,
from the properties of $R_i$ and~$S_i^*$.
\end{remark}

In addition to $R_i$, $R_i^*$, $S_i$ and
$S_i^*$ we also consider the operators
$M_i$ and $M_i^*$, defined on Hermite polynomials by
$$
M_i=x_i\, M(\cL) \qquad\qquad M^*_i=M(\cL)\, x_i.
$$

\begin{corollary}
The operators $M_i$ and $M_i^*$ extend to bounded operators on $\lp{\ga}$
for all $p$ in $(1, \infty)$ in any dimension. 
For $d=1$, they are bounded from $\hu{\ga}$ to $\lu{\ga}$ 
and from $\ly{\ga}$ to  $\BMO{\ga}$.  For $d\geq 2$, they 
are unbounded both from $\hu{\ga}$ to $\lu{\ga}$ 
and from $\ly{\ga}$ to  $\BMO{\ga}$. 
\end{corollary}

\begin{proof}
Since $M_i =(1/2) (R_i+S^*_i)$ and $M^*_i=(1/2) (R^*_i+S_i)$,
the $\lp{\ga}$ boundedness of $M_i$ and $M_i^*$
follows from that of $R_i$, $S_i$, $R_i^*$ and $S_i^*$
mentioned in the introduction.  The remaining parts
of the corollary 
are straightforward consequences of Theorem~\ref{t: main}.
\end{proof}

The proof of Theorem~\ref{t: main} occupies the remaining part of
the paper.
We first deal with $R_i$ and $S_i^*$. In Section~\ref{R}, 
it is proved that in dimension one, 
$R_i$ is bounded from $\hu{\ga}$ to $\lu{\ga}$ and $S_i^*$ 
from $\ly{\ga}$ to $BMO(\ga)$.  
Section~\ref{S} deals with the positive
statements concerning these two operators for $d\ge 2$, as indicated in the
table. The negative results claimed for  $R_i$ and $S_i^*$ are
proved in Section~\ref{U}. Finally, a rather simple duality argument
in Section~\ref{D}
gives all the results for the remaining operators $R_i^*$ and $S_i.$

\section{Boundedness of $R_1$ and $S_1^*$ in one dimension}  \label{R}

Let $d=1$.  We first prove that $R_1$ is bounded
from $\hu{\ga}$ to $\lu{\ga}$.
Because of (\ref{f: HIC H1L1}) and (\ref{eq:kern}),
it suffices to show that there exists a constant $C$
such that for all balls $B$ in~$\cB_1$ and all points $y$ in $B$
\begin{equation}\label{basic2}
r_B \,  \int_{(2B)^c} \!\!\! \wrt\la(x)\, 
\Bigmod{\int_0^1 \frac{r\,\partial_y\bigl(\psi\, 
\e^{-\phi^2}\bigr)}{1-r^2} \, \wrt {\rho(r)} }
\leq C.
\end{equation}
Next, observe that
\begin{align*}\notag
\partial_y\bigl(\psi\ \e^{-\phi^2}\bigr) 
& = -\frac{1}{\sqrt{1-r^2}}\,\e^{-\phi^2}-2\phi\,
    \psi \frac{r}{\sqrt{1-r^2}}\,\e^{-\phi^2} \\
& = F_1(r,x,y)+F_2(r,x,y).
\end{align*}
To deal with $F_2$, we first note that 
$\partial \phi/\partial r = - \psi/(1-r^2)$, so that
$$
F_2(r,x,y) = -r\, \sqrt{1-r^2}\, \partial \e^{-\phi^2}/\partial r.
$$ 
Integrating by parts, we get 
\[
\int_{0}^1 \frac{(-\log r)^{-1/2}}{1-r^2}\,F_2(r,x,y)\,\wrt r =
\int_{0}^1 \frac{\partial}{\partial r}
\frac{r(-\log r)^{-1/2}}{\sqrt{1-r^2}}\,\e^{-\phi^2}\,\wrt r.
\] 
The integral to the right here is easily seen to be bounded by
\[
C\, \int_{0}^1 (1-r)^{-2}\, \e^{-\phi^2} \wrt r.
\]
Considering also $F_1$, we can then estimate the expression in 
(\ref{basic2}) by
\begin{equation}
  \label{eq:upper}
C\, r_B \int_{(2B)^c} \int_{0}^1 (1-r)^{-2}\, \e^{-\phi^2}\,\wrt r\wrt\la(x).
\end{equation}
Here we first integrate in $r$ only over  $1-r_B/(2(1+|y|)) < r < 1$,
which implies
\[
|ry-x| = |(r-1)y+y-x| \geq |y-x| - \frac{r_B|y|}{2(1+|y|)} \geq
|y-x| - \frac{r_B}2 \geq \frac{|y-x|}2
\]
and thus
$|\phi| \geq |y-x|/(2\sqrt{1-r})$.
This part of (\ref{eq:upper})
can thus be estimated by the following expression, where
we make the change of variables  $t = (1-r)/|x-y|^2$,

\begin{align*}
C\, r_B \int_{(2B)^c} \int_{1-r_B/(2(1+|y|))}^1  
\!\! &\frac{\e^{-c|x-y|^2/(1-r)}\,}{(1-r)^2}  \wrt r \wrt\la(x) \\
&\leq r_B  \int_{(2B)^c} \frac{\wrt\la(x)}{|x-y|^2} 
\int_0^\infty t^{-2}\e^{-c/t} \wrt t \leq C.
\end{align*}
In the remaining part of  (\ref{eq:upper}), we switch
the order of integration and find
\begin{align*}
r_B \int_{(2B)^c}\!\!\! \wrt\la(x)\, &\int_{0}^{1-r_B/(2(1+|y|))}   
     \frac{\e^{-\phi(r,x,y)^2}}{(1-r)^2}  \wrt r     \\  
&\leq r_B \int_{0}^{1-r_B/(2(1+|y|))} \!\!\frac{\wrt r}{(1-r)^2} 
    \int_{\BR^d} \e^{-\phi(r,x,y)^2} \wrt\la(x) \\
&\leq C \, r_B \int_{0}^{1-r_B/(2(1+|y|))}\frac{\wrt r}{(1-r)^{3/2}},
\end{align*}

which is easily seen to be bounded by a constant independent
of $y$ in $B$ and of $B$ in $\cB_1$ .  The proof of the
boundedness of $R_1$ from $\hu{\ga}$ to $\lu{\ga}$ is complete.

\medskip
Now we prove that if $d=1$, then $S_1^*$ is bounded from $\ly{\ga}$
to $\BMO{\ga}$.
Recall that the functions $\phi$ and $\psi$ are now scalar valued.
By (\ref{f: HIC LIBMO}) and (\ref{eq:kerntilde}) it suffices to show that
there exists a constant $C$, independent of $x$ in $B$ and of 
$B$ in $\cB_1$, such that 
\begin{equation} \label{equiv}
r_B \, 
\int_{(2B)^c} \!\! \wrt \la(y)\, \Bigmod{\int_0^1 \frac{1}{1-r^2} \, 
\partial_x \bigl(\phi\, \e^{-\psi^2} 
+{(1-r^2)} \, x \e^{-y^2}\bigr) \wrt \rho(r)}  
\leq C.
\end{equation}
This  is analogous to  (\ref{basic2}), with the following
three differences. The variables $x$ and $y$, and thus also $\phi$
and $\psi$, are swapped; there is an extra factor
$1/r$; there is an extra term ${(1-r^2)} x \e^{-y^2}$.

For that part of the expression in (\ref{equiv}) obtained by integrating
only over $1/2 < r < 1$, we can copy the argument used to verify 
(\ref{basic2}) above. The extra factor $1/r$ is harmless
in this interval, and so is the term ${(1-r^2)} x \e^{-y^2}$. 
In the integration by parts with respect to $r$, one will now get
an integrated term at $r=1/2$, equal to $C\,\e^{-\psi(1/2,x,y)^2}$ and
easy to handle.

Consider then the integral over $0 < r < 1/2$. We have
\begin{equation} \label{eq:derivative}
\partial_x \bigl(\phi \,\e^{-\psi^2} +{(1-r^2)} \,x \,\e^{-y^2}\bigr)
= -(1-r^2)\, \Bigl[\frac{\e^{-\psi^2} }{(1-r^2)^{3/2}} - \e^{-y^2}\Bigr]
-\frac{2\phi\psi r}{(1-r^2)^{1/2}} \, \e^{-\psi^2}.
 \end{equation}
Notice first that $ry-x = r(y-rx) + (r^2-1)x$, so that
\begin{equation} \label{eq:little}
 |\phi| 
\le |\psi| + |x|.
\end{equation}
This allows us to estimate the last term in (\ref{eq:derivative}) by 
 $C\,r(1+|x|)\e^{-\psi^2/2}$ for $0 < r < 1/2$. 
The corresponding part of the expression in 
(\ref{equiv}) can then be seen not to be larger than
$C\, r_B (1+|x|) \le C$, since $x \in B$.
For the first term in the right-hand side of (\ref{eq:derivative}), we write
\[
\frac{  \e^{-\psi^2} }{(1-r^2)^{3/2}} - \e^{-y^2}
= \frac{ \e^{-\psi^2} -  \e^{-y^2} }{(1-r^2)^{3/2}}
    +\Bigl[\frac{ 1 }{(1-r^2)^{3/2}} - 1\Bigr] \,  \e^{-y^2}.
\]
The last  term here is easily seen to be harmless and will be neglected. 
For the first term to the right, we observe that 
$\e^{-y^2}  = \e^{-\psi(0,x,y)^2}$, so that
\[
\e^{-\psi(r,x,y)^2} -  \e^{-y^2} = \int_0^r \frac{\partial}{\partial s}
\e^{-\psi(s,x,y)^2}\wrt s = \int_0^r 
\frac{2 \phi(s,x,y) \psi(s,x,y)}{1-s^2}\e^{-\psi(s,x,y)^2}\wrt s.
\]
Since $s<r<1/2$
here,  (\ref{eq:little}) implies that the last integral is at most
 \begin{equation*}
   \label{eq:s-integral}
 C\, \int_0^r (1 + |x|)\, \e^{-\psi(s,x,y)^2/2}\wrt s.
 \end{equation*} 
Integrating this expression with respect to $y$, we get at most
$C\,r(1+|x|)$. This will give a contribution to the expression in
(\ref{equiv}) which is at most 
$C\,r_B (1+|x|) \int_0^{1/2} (- \log r)^{-1/2} \wrt r$, and this is 
uniformly bounded with respect to $x$ in $B$ and $B$ in $\cB_1$.

This concludes the proof that $S_1^*$ is bounded from 
$\ly{\ga}$ to $\BMO{\ga}$. 

\section{Boundedness of $R_i$ and $S_i^*$ in higher dimensions}
\label{S}

The boundedness of $R_i$ from $\ly{\ga}$ to  $\BMO{\ga}$ for any $d$
was proved in \cite[Theorem 7.2~\rmii]{MM}. 

To prove that $S_i^*$ is bounded from $\hu{\ga}$ to $\lu{\ga}$ in
all dimensions,
we follow the lines of the proof of \cite[Theorem 7.2]{MM}.
By (\ref{f: HIC H1L1}) and (\ref{eq:kerntilde}),
it is enough to verify that there exists a constant $C$, independent
of $y$ in $B$ and of $B$ in $\cB_1$ such that 
\begin{equation}\label{basic3}
r_B \, \int_{(2B)^c} \Bigmod{\int_0^1 \nabla_y
\Bigl(\frac{\phi_i \, \e^{-|\phi|^2}}{(1-r^2)^{(d+1)/2}}
+ x_i \, \e^{-|x|^2}\Bigr) \wrt \rho(r)} \wrt\la(x)
\leq C.
\end{equation}
It is straightforward to verify that all the components of the gradient 
with respect to $y$ that appear here are bounded in absolute value by
\[
C\, \frac{r\, (1+|\phi|^2)}{(1-r^2)^{(d+2)/2}} \,  \e^{-|\phi|^2} 
\le C\,\frac{r}{(1-r^2)^{(d+2)/2}} \, \e^{-|\phi|^2/2}.
\]
Thus the left-hand side of (\ref{basic3}) is no larger than
\begin{equation} \label{eq:doubleint}
C\, r_B
\int_{(2B)^c} \wrt\la(x) \int_0^1   
\frac{r\, \e^{-|\phi|^2/2}}{(1-r^2)^{(d+2)/2}}  \wrt \rho(r).
\end{equation}
We split the inner integral into integrals over
$(0, r_0)$ and $(r_0, 1)$, with  $r_0 = 1 - r_B^2/4> 1/2$.
In each of the two resulting double integrals,
we switch the order of integration. That part of (\ref{eq:doubleint})
corresponding to $(0, r_0)$ will be less than
\begin{align*}
\quad C\, r_B \int_0^{r_0}  \frac{(-\log r)^{-1/2}}{(1-r^2)^{(d+2)/2}}  
      \wrt r \, \int_{\BR^d} \e^{-|\phi|^2/2}
      & \wrt\la(x)  \\
& \le  C\, r_B  \int_0^{r_0} (-\log r)^{-1/2} (1-r^2)^{-1}  \wrt r  \\  
& \le  C \, r_B \, (1-r_0)^{-1/2}   \\ 
& \le  C.
\end{align*}  
To deal with the integral over $(r_0, 1)$, 
observe that $y \in B$ implies $|y| \le c_B + r_B
\le 2 r_B^{-1}$, in view of the definition of $\cB_1$. For $x \notin 2B$
and  $r_0 < r < 1$ we then have 
\[
|ry-x| = |y-x-(1-r)y| > r_B -|y|r_B^2/4 \ge r_B/2.
\]
The remaining part of (\ref{eq:doubleint})  is thus at most
\begin{equation} \label{eq:switch}
 C\, r_B \int_{r_0}^1 \frac{r}{(1-r^2)^{(d+2)/2}} \wrt \rho(r) \,
\int_{|ry-x|>r_B/2} \e^{-|\phi|^2/2} \wrt\la(x). 
\end{equation}
The inner integral in  (\ref{eq:switch}) is dominated by
\[
(1-r^2)^{d/2}\, \int_{|z|>r_B/(2\sqrt{1-r^2})} \e^{-|z|^2} \wrt z
< C\, (1-r^2)^{d/2}\,\left(r_B/\sqrt{1-r^2}\right)^{-2},
\]
and the whole expression (\ref{eq:switch})  by 
\[
C\, r_B^{-1}  \int_{r_0}^1 (1-r)^{-1/2}\wrt r \le C\, r_B^{-1} (1-r_0)^{1/2} \le C.
\]
This concludes the proof of the 
boundedness of $S_i^*$ from $\hu{\ga}$ to $\lu{\ga}$.

\section{Unboundedness results}
\label{U}

In this section, we shall prove that if $d\geq 2$, then 
$R_1$ is unbounded from 
$\hu{\ga}$ to $\lu{\ga}$ and $S_1^*$ is unbounded from 
$\ly{\ga}$ to $\BMO{\ga}$.
The proofs for $R_2,\ldots,R_d$ and $S_2^*, \ldots, S_d^*$
are similar and are omitted. 

Consider first $R_1$.
Let $Q$ be the cube with centre $(\xi,0,...,0)$ and 
sidelength $2\xi^{-1}$, for some large $\xi>0$. 
Denote by $Q_+$ and $Q_-$ those halves  of $Q$
where $y_2>0$ and  $y_2<0$, respectively. 
{Set $b=\ga(Q)^{-1}(\1_{Q_+}-\1_{Q_-})$.  
By Remark~\ref{rem: eq norm BMO}, $\norm{b}{H^1(\gamma)}$ is bounded
independently of $\xi$.}    
Writing $\tilde y=(y_1,-y_2,y_3,\ldots,y_d)$, one can easily see from 
(\ref{eq:kern}) that
\begin{align}\label{Ra}
\e^{-\mod{x}^2}\, R_1 a(x)
& =  \frac2{\sqrt{\pi}\,\ga(Q)}\int_{Q_+} \int_{0}^1
       \frac{r\, (y_1-rx_1)}{(1-r^2)^{(d+2)/2}}\, 
       \big(\e^{-|\phi(r, x,  y)|^2}-\e^{-|\phi(r, x, \tilde  y)|^2}\big)
       \wrt \rho(r) \wrt\ga(y)\notag \\
& =  \frac2{\sqrt{\pi}\,\ga(Q)}\int_{Q_+} \int_{0}^1
       \frac{r\, (y_1-rx_1)}{(1-r^2)^{(d+2)/2}} \,
       \e^{-|\phi(r, x, y)|^2}\, \big(1-\e^{-\tau(r,x_2,y_2)}\big)
       \wrt \rho(r) \wrt\ga(y);
\end{align}
here and in the rest of the proof we write 
$\tau(r,x_2,y_2)$ in place of $4rx_2y_2/(1-r^2)$.
Similarly, we write $\eta(\xi,x_1)$ instead of $\sqrt{(\xi-x_1)/\xi}$.
We shall evaluate $R_1a$ at points $x$ in the set
$$
A_Q =
\Bigl\{x\in\BR^d: x_1\in \Bigl(\xi - 1, \xi-\frac4\xi\Bigr), \, 
x_2 \in \Bigl[\frac{\eta(\xi,x_1)}{2} , \eta(\xi,x_1)\Bigr],
|x_k| \leq \eta(\xi,x_1),  k=3,...,d\Bigl\}.
$$
Observe that if $x\in A_Q$, $y\in Q_+$ and $r\in [0,1]$, then $y_1-rx_1>0$
and the integrands  in (\ref{Ra}) are positive.
With $\xi$ large and $x \in A_Q$, the integration in $r$ will be 
restricted to the interval 
\begin{equation} \label{f: Int}
I(x_1) = \Bigset{r\in(0,1): \Bigmod{r-\frac{x_1}{\xi}}
<\frac{\eta(\xi,x_1)}{2\xi} }\subset (1/2, 1).
\end{equation}
For  $x \in A_Q$ a lower bound for $\e^{-\mod{x}^2}\,R_1a(x)$ is
$$
\frac{c}{\ga(Q)}\int_{Q_+} \int_{I(x_1)} \frac{y_1-rx_1}{(1-r^2)^{(d+3)/2}}\,
\e^{-|\phi|^2}\ \big(1-\e^{-\tau(r,x_2,y_2)}\big) \wrt r\wrt\ga(y). 
$$
Therefore, 
\begin{multline}\label{difference}
\norm{R_1a}{1}
\ge \int_{A_Q} R_1a(x) \, \e^{-\mod{x}^2}\wrt \lambda(x) \\  
\ge \frac{c}{\ga(Q)}\int_{A_Q}\int_{Q_+} \int_{I(x_1)} 
\frac{y_1-rx_1}{(1-r^2)^{(d+3)/2}}\, \e^{-|\phi|^2}\, 
\big(1-\e^{-\tau(r,x_2,y_2)}\big)\wrt r\wrt\ga(y) \,\wrt\lambda(x). 
\end{multline}
To prove that the operator $R_1$ is unbounded from $\hu{\ga}$ to $\lu{\ga}$, 
we only need to show that the last expression here
tends to infinity with $\xi$. 
We claim that if $\xi$ is sufficiently large, 
one has for $x$ in $A_Q$, $y$ in $Q_+$ and  $r$ in $I(x_1)$ 
\begin{align}
 1-r^2 
& \sim \eta(\xi,x_1)^2   \label{dis1} \\ 
y_1-rx_1
&\ge c\  (\xi-x_1) \label{dis2}\\
\e^{-|\phi(r, x, y)|^2}
&\ge c     \label{dis3}\\
1-\e^{-\tau(r,x_2,y_2)}
&\sim y_2/\eta(\xi,x_1).  \label{dis4}    
\end{align}
Deferring momentarily their proofs, 
we show how these inequalities yield the desired conclusion. 
Indeed, applying them  and observing that 
$$
\frac{1}{\ga(Q)}\int_{Q_+}y_2\wrt\ga(y) \ge\, \frac{c}\xi,
$$
we see from (\ref{difference}) that
\begin{align*}
\norm{R_1a}{1}
& \ge  \frac{c}{\ga(Q)}\int_{A_Q}\int_{Q_+}
        \eta(\xi,x_1)^{-(d+2)}\, \xi \, y_2\, \mod{I(x_1)}  
        \wrt\ga(y)\,\wrt\lambda(x) \\ 
& \ge  c\,\int_{A_Q}\eta(\xi,x_1)^{-(d+1)}\, \xi^{-1} \wrt\lambda(x) \\ 
& =    c\,\int_{\xi - 1}^{\xi-4/\xi} (\xi-x_1)^{-1} \wrt x_1\\
& \ge  c\,\log \xi,
\end{align*}
which tends to infinity with $\xi$ as desired.\par
It remains  to prove the inequalities (\ref{dis1})--(\ref{dis4}). 
Observe first that for $x \in A_Q$ one has $4/\xi < \xi - x_1$.
Considering also the geometric mean of these two quantities, we have
\begin{equation}
  \label{eq:star}
  4/\xi < 2\, \eta(\xi,x_1) <\xi - x_1.
\end{equation}

To verify (\ref{dis1}), write for $r \in I(x_1)$
$$
1-r^2
\sim 1-r 
=\eta(\xi,x_1)^2 -\Bigl(r-\frac{x_1}\xi\Bigr).
$$
Observe that $|r-x_1/\xi| < \eta(\xi,x_1)^2/4$, because of
the definition of $I(x_1)$ and (\ref{eq:star}), and  (\ref{dis1})
follows.

For (\ref{dis2}), notice  that similarly
$$
y_1-rx_1 > \xi - \frac1\xi
- \Bigl(\frac{x_1}\xi +  \frac{\xi-x_1}{4\xi}\Bigr) x_1
> \xi - \frac{\xi-x_1}4 -\frac{x_1^2}{\xi} - \frac{\xi-x_1}4
> \frac{\xi^2-x_1^2}{\xi} - \frac{\xi-x_1}2
>  \frac{\xi-x_1}2. 
$$

Aiming at (\ref{dis3}), we write
\[
\mod{x_1-ry_1} \le \mod{r(y_1-\xi)} +  \mod{\left(r - \frac{x_1}\xi\right)\xi}
<  \frac1\xi + \frac{1}{2\xi}\eta(\xi,x_1)\,\xi <
\eta(\xi,x_1) \le C\,\sqrt{1-r^2},
\]
the last two inequalities in view of (\ref{eq:star}) and (\ref{dis1}),
respectively. Since $x\in A_Q$ and $y\in Q_+$, we similarly get
 for $k=2,... ,d$
$$
\mod{x_k-ry_k}
\le \mod{x_k}+\mod{y_k}
\le \eta(\xi,x_1) +\xi^{-1}
\le C\,\sqrt{1-r^2}.
$$
 Thus $\mod{x-ry}^2/(1-r^2)\le C$, which implies (\ref{dis3}).

Finally, we prove (\ref{dis4}).  Because of (\ref{dis1}) and the facts that
$x\in A_Q$ and $y\in Q_+$, one has
$$
\tau(r,x_2,y_2) \sim \frac{\xi}{\xi-x_1}\, \eta(\xi,x_1)\, y_2
= \eta(\xi,x_1)\, y_2 < \frac1{\sqrt{\xi(\xi-x_1)}} < 1,
$$
and (\ref{dis4}) follows.

This concludes the proof that $R_1$ is unbounded
from $\hu{\ga}$ to $\lu{\ga}$.

\medskip
Next, we prove that in higher dimension, $S_1^*$ is unbounded from $\ly{\ga}$
to $\BMO{\ga}$.
We shall modify slightly the counterexample 
constructed above.
The kernel is
\[
k_{S_1^*}(x,y)= \frac2{\sqrt{\pi}}\, \e^{|x|^2} \int_0^1
\Bigl[\frac{x_1 - ry_1}{(1-r^2)^{(d+2)/2}}  \,
\e^{-|\phi|^2} + x_1 \e^{-|x|^2}\Bigr] \wrt \rho(r).
\]
With $\xi > 0$ large, we let the cube $Q$, its subset  $Q_+$ and the region 
$A_Q$ be as above.  But now  $f$ will be the characteristic
function of $A_Q$, and $S_1^* f$ will be evaluated in $Q$. 
For $x \in Q_+$ we write $\tilde{x} = (x_1, -x_2, x_3, ..., x_d)$.
We shall verify that
\begin{equation}
  \label{nonbmo}
\frac1{\ga(Q)} \int_{Q_+} (S_1^* f({x}) - S_1^* f(\tilde x))
\wrt \ga(x)
 \end{equation}
tends to $+\infty$ as $\xi$ tends to $+\infty$. This would clearly mean that 
$\norm{S_1^* f}{BMO}$ also tends to infinity, although
$\norm{f}\infty = 1$ for all $\xi$. Notice that for $x$ in $Q_+$
\[
S_1^* f({x}) - S_1^* f(\tilde x) = \frac2{\sqrt{\pi}}\, \int_{A_Q} \int_0^1
\frac{x_1-ry_1}{(1-r^2)^{(d+3)/2}} 
\e^{-\mod{\psi}^2}\, \big(1-\e^{-\tau(r,x_2,y_2)}\big) 
\wrt \rho(r) \wrt\la(y);
\]
recall that we write $\tau(r,x_2,y_2)$ in place of $4rx_2y_2/(1-r^2)$.
Here the integrand is positive. We now integrate  with respect to
$\wrt \ga(x)$ over $x \in Q_+$, and restrict the integral in $r$
to the interval $I(y_1) \subset (1/2, 1)$, defined like $I(x_1)$ 
(see (\ref{f: Int})).
For such $r$, one can replace the measure $\wrt \rho(r)$ by
$(1-r^2)^{-1/2} \wrt r$.
 This leads to a lower estimate of the quantity 
(\ref{nonbmo}) in terms of  a triple integral rather similar to that
in (\ref{difference}). The only differences are that $x$ and $y$ 
have switched roles and that we now have also an insignificant factor $1/r$. 
Thus the integral will tend to
$\infty$ with $\xi$, like that in (\ref{difference}). 

This concludes the proof that $S_1^*$ is unbounded from $\ly{\ga}$
to $\BMO{\ga}$.

\section{Duality arguments}
\label{D}

All the statements contained in Theorem~\ref{t: main}
that have not been proved yet will follow from those
proved above by a simple duality argument,
based on the following elementary lemma. 
We point out that by a result of \cite{MMS2},  
$\ld{\ga}$ is a subspace of $\hu{\ga}$, 
which is dense because it contains all the atoms. 

\begin{lemma}  \label{l: duality}
Suppose that $T$ is a bounded operator on $\ld{\ga}$.  Then 
$T$ extends to a bounded operator from $\hu{\ga}$ to $\lu{\ga}$
if and only if its (Hilbert space) adjoint operator $T^*$
is bounded from $\ly{\ga}$ to $\BMO{\ga}$.
\end{lemma} 

\begin{proof}
Suppose that $T$ extends to a bounded operator from $\hu{\ga}$
to $\lu{\ga}$. 
Then for $g$ in~$\ly{\ga}$
\begin{equation} \label{f: norm T*f in BMO}
\norm{T^*g}{BMO} 
 = \sup\,\{ \mod{\prodo{f}{T^*g}}: f \in \ld{\ga}, \ 
       \norm{f}{H^1}\le 1\},
\end{equation}
and 
$$
\begin{aligned}
\mod{\prodo{f}{T^*g}}
& =   \mod{\prodo{Tf}{g}} \\
& \le \norm{Tf}{1}\, \norm{g}{\infty}  \\
& \le \norm{T}{H^1,L^1} \, \norm{f}{H^1}\, \norm{g}{\infty}.
\end{aligned}
$$
Thus $T^*$ is bounded from $\ly{\ga}$ to $BMO(\ga)$.

The converse implication is proved similarly.  We omit the details.
\end{proof}

Now, the boundedness in one dimension of $R_1^*$ from $\ly{\ga}$ to $\BMO{\ga}$ 
follows  by Lemma~\ref{l: duality}
from that of $R_1$ from $\hu{\ga}$ to $\lu{\ga}$, proved in 
Section~\ref{R}.   Similarly, the
unboundedness in dimension $d\geq 2$ of $R_i^*$ from $\ly{\ga}$ to $\BMO{\ga}$
follows from the unboundedness of $R_i$ from $\hu{\ga}$ to $\lu{\ga}$
proved in Section~\ref{U}, 
and the boundedness of $R_i^*$ from $\hu{\ga}$ to  $\lu{\ga}$ for any $d$
follows from that of $R_i$ from $\ly{\ga}$ to $\BMO{\ga}$
proved in Section~\ref{S}.

A similar argument gives the boundedness of $S_i$ 
from $\hu{\ga}$ to $\lu{\ga}$ when $d=1$,
that of $S_i$ from $\ly{\ga}$ to  $\BMO{\ga}$ for any $d$,
and the unboundednes of  $S_i$ from $\hu{\ga}$ to $\lu{\ga}$ 
when $d\geq 2$.

The theorem is completely proved.

\end{document}